\documentclass{amsart}

\usepackage{cite}
\usepackage{graphicx}
\usepackage{amsthm,amscd,amsmath,amssymb,amsfonts,enumerate}
\usepackage{braket}

\newtheorem{theorem}{Theorem}

\title[Representation Theory, Voting Theory, and Game Theory]{Representation Theory of the Symmetric Group \\ in Voting Theory and Game Theory}

\author{Karl-Dieter Crisman}
\address{Department of Mathematics and Computer Science, Gordon College, 255 Grapevine Road, Wenham, Massachusetts 01984}
\curraddr{}
\email{karl.crisman@gordon.edu}
\thanks{}

\author{Michael E.~Orrison}
\address{Department of Mathematics, Harvey Mudd College, 301 Platt Boulevard, Claremont, California 91711}
\curraddr{}
\email{orrison@hmc.edu}
\thanks{}

\date{August 7, 2015}

\begin{document}
\maketitle


\begin{abstract}
This paper is a survey of some of the ways in which the representation theory of the symmetric group has been used in voting theory and game theory.  In particular, we use permutation representations that arise from the action of the symmetric group on tabloids to describe, for example, a surprising relationship between the Borda count and Kemeny rule in voting.  We also explain a powerful representation-theoretic approach to working with linear symmetric solution concepts in cooperative game theory.  Along the way, we discuss new research questions that arise within and because of the representation-theoretic framework we are using.  
\end{abstract}


\section{Introduction}

Symmetry arises in many ways in voting theory and game theory. In this paper, we highlight one such appearance by surveying some of the ways in which the \emph{representation theory} of the symmetric group has been used in these fields.  Our primary goal is to show how certain well-understood permutation representations of the symmetric group can be used to make sense of foundational ideas in voting theory and game theory, and how using these representations can in turn help researchers formulate novel questions and meaningful generalizations.

Among other examples, we use the representation theory of the symmetric group to describe a surprising relationship between the Borda count and the Kemeny rule in voting.  We also use it to construct and make sense of certain infinite families of solution concepts in cooperative game theory. Along the way, we describe much of the representation-theoretic framework used in voting theory and game theory, and we discuss several new research questions that arise therein.

The permutation representations on which we focus are those associated to the usual action of the symmetric group $S_n$ on simple combinatorial objects called tabloids.  Representation theory experts will have no trouble seeing that we only scratch the surface on what could be done with these representations in voting theory and game theory.  Another goal for this paper is therefore to encourage readers to contribute to these fields by going well beyond what we present here.

As a brief introduction to the role that symmetry can play in voting, consider a situation in which voters are asked to choose their favorite candidate $A$, $B$, or $C$.  If we end with a profile of $(13,4,7)$ for $A,B,C$ respectively, then it is natural to say this is essentially the same profile as $(9,0,3)$, because there are four votes for each candidate that ``cancel" each other.  We can go further, though, by writing $(13, 4, 7) = (8,8,8) + (5,-4,-1)$ and noting that the first vector on the right captures how many voters there were, while the second vector captures the different levels of support the candidates received, which in the end is what really matters.  

Interestingly, such decompositions profiles arise by having the symmetric group $S_3$ act on the candidates by permuting their labels, which in turn induces an action on the set of profiles.  Furthermore, by having the symmetric group act on the labels of candidates in much more complicated voting situations, similar decompositions of profiles arise, and we are able to use what we know about the resulting subspaces of voting data to say something worthwhile (as we will demonstrate in Section 3) about how profiles are used by different kinds of voting procedures.  

In the next section, we introduce most of the notation we will use for the rest of the paper, and we motivate the use of tabloids for indexing the kind of data with which we will be working.  Although we do not assume our readers will be familiar with the ideas we will be presenting from voting theory and game theory, we will assume that our readers have a basic working knowledge of the representation theory of finite groups (see, for example, \cite{JamesLiebeck, Serre}).

Finally, we encourage interested readers to view the ideas presented in this paper within the larger framework of harmonic analysis on finite groups \cite{ClausenBaum, CSSTHarm, Diaconis, Terras}.  In our opinion, doing so makes it much easier to see how the representation theory of the symmetric group (and other finite groups) has been  applied outside of voting theory and game theory, in fields, for example, such as statistics \cite{Diaconis, DiaconisEriksson, Marden, BargOrrLinearRank} and machine learning \cite{HuangETALProbInf, KondorQAP, KondorMulti}.  Doing so will also reveal that the approach we are taking in this paper has already found success in a variety of other settings.


\section{Background}

In this section, we define tabloids and their associated permutation representations of the symmetric group.  We also describe how tabloids can be used to index the kind of voting and game-theoretic data with which we will be working.  Good references for the material in this section are \cite{JamesKerber} and \cite{Sagan}.

Let $n$ be a positive integer.  A \emph{composition} of $n$ is a list $\lambda = (\lambda_1, \dots, \lambda_m)$ of positive integers whose sum is $n$.  If it is also the case that $\lambda_1 \ge \cdots \ge \lambda_m$, then we say $\lambda$ is a \emph{partition} of $n$.  For example, $(2, 3, 1, 3)$ is a composition of $9$, and $(4,2,2,1)$ is a partition of $9$.  If $\lambda$ is a composition, then we will denote by $\overline{\lambda}$ the partition obtained by reordering the numbers in $\lambda$ so that they form a non-increasing list.  For example, if $\lambda = (2, 3, 1, 3)$, then $\overline{\lambda} = (3,3,2,1)$.  

\begin{figure}[h]
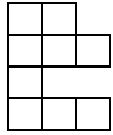

\begin{center}
\begin{tabular}{| c | c | c |} \cline{1-2}
\hspace{.04in} & \hspace{.04in} & \multicolumn{1}{| c}{ } \\ \cline{1-3}  
\  &   & \hspace{.04in}  \\ \cline{1-3}
\  & \multicolumn{2}{| c}{} \\ \cline{1-3}
\  &   &   \\ \cline{1-3} 
\end{tabular}
\end{center}
\caption{The Young diagram of shape $(2,3,1,3)$.}
\label{youngdiagram}
\end{figure}

If $\lambda$ is a composition of $n$, then the \emph{Young diagram of shape $\lambda$} is the left-justified array of boxes that has $\lambda_i$ boxes in its $i$th row (see Figure \ref{youngdiagram}).  If we fill these boxes with the numbers $1, \dots, n$ without repetition, then we create a \emph{Young tableau of shape $\lambda$}.  Two tableaux of shape $\lambda = (\lambda_1, \dots, \lambda_m)$ are then said to be \emph{row equivalent} if they have the same set of $\lambda_1$ numbers in the first row, the same set of $\lambda_2$ numbers in the second row, and so on.  An equivalence class of tableaux is then called a \emph{tabloid of shape $\lambda$}.

\begin{figure}[h]
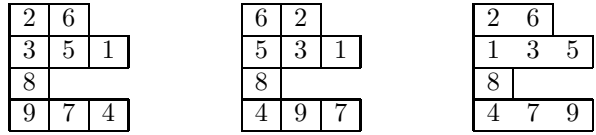

\begin{center}
\begin{tabular}{| c | c | c |} \cline{1-2}
2 & 6 & \multicolumn{1}{| c}{ } \\ \cline{1-3}  
3 & 5 & 1 \\ \cline{1-3}
8 & \multicolumn{2}{| c}{} \\ \cline{1-3}
9 & 7 & 4 \\ \cline{1-3} 
\end{tabular}
\hspace{.5in}
\begin{tabular}{| c | c | c |} \cline{1-2}
6 & 2 & \multicolumn{1}{| c}{ } \\ \cline{1-3}  
5 & 3 & 1 \\ \cline{1-3}
8 & \multicolumn{2}{| c}{} \\ \cline{1-3}
4 & 9 & 7 \\ \cline{1-3} 
\end{tabular}
\hspace{.5in}
\begin{tabular}{| c c c |} \cline{1-2}
2 & 6 & \multicolumn{1}{| c}{ } \\ \cline{1-3}  
1 & 3 & 5 \\ \cline{1-3}
8 & \multicolumn{2}{| c}{} \\ \cline{1-3}
4 & 7 & 9 \\ \cline{1-3} 
\end{tabular}
\end{center}
\caption{Two row equivalent tableaux and their tabloid.}
\label{tabloid}
\end{figure}

We often denote a tabloid by first forming a representative tableau and then removing the vertical dividers within each row (see Figure \ref{tabloid}).  For convenience, we will usually choose the representative tableau whose entries in each row are in ascending order, and we will read the entries from top to bottom when putting such representatives and their associated tabloids in lexicographic order.  

Suppose $\lambda$ is a composition of $n$.  Let $X^\lambda$ denote the set of tabloids of shape $\lambda$, and let $M^\lambda$ denote the vector space of real-valued functions defined on $X^\lambda$.  If $x \in X^\lambda$, then $\mathbf{f}_x$ will denote its associated indicator function with the property that $\mathbf{f}_x(x) = 1$ and $\mathbf{f}_x(y) = 0$ for all $y \ne x$.  Note that the set of such indicator functions forms an orthonormal basis for $M^\lambda$ with the respect to the usual inner product $\langle \cdot , \cdot \rangle$ on $M^\lambda$, which is defined by setting
\[
\langle \mathbf{f}, \mathbf{g} \rangle = \sum_{x \in X^\lambda} \mathbf{f}(x) \mathbf{g}(x)
\]
for all $\mathbf{f}, \mathbf{g} \in M^\lambda$.  We will call this basis of indicator functions the \emph{usual basis} of $M^\lambda$, and when necessary we will assume it has been ordered with respect to the lexicographic ordering of the associated tabloids in $X^\lambda$.  

In this paper, we will be interested in a variety of real-valued functions defined on $X^\lambda$ for different choices of $\lambda$.  We will also be interested in certain linear transformations defined on $M^\lambda$.  

For example, consider an election involving $n$ candidates labeled $1$ through $n$.  If we ask the voters to tell us their top favorite $\lambda_1$ candidates, then their next top favorite $\lambda_2$ candidates, and so on, then we can interpret each of their responses as a tabloid $x \in X^\lambda$.  For example, if $\lambda = (2,3,1,3)$, and a voter chooses the tabloid in Figure \ref{tabloid}, then we will take that to mean candidates 2 and 6 are her top two favorite candidates, but we do not know if the voter prefers one of these two candidates over the other.  The function $\mathbf{f} \in M^\lambda$ in which we are then interested is simply the one for which $\mathbf{f}(x)$ is the total number of voters who chose tabloid $x \in X^\lambda$.

We will mostly focus on three kinds of compositions in this paper.  The first composition is the all-ones partition $(1,\dots,1)$.  The tabloids in $X^{(1, \dots, 1)}$ naturally correspond to all of the possible permutations of the numbers $1,\dots, n$.  Therefore, $M^{(1,\dots,1)}$ is helpful when describing voting situations in which voters are asked to provide a full ranking of candidates.  The second kind of composition has the form $(k, n-k)$ where $1 \le k \le n-1$.  In this case, we will use the tabloids in $X^{(k,n-k)}$ to index  $k$-element subsets of $\{ 1, \dots, n \}$, where each $k$-element subset corresponds to the tabloid whose top row contains all of its elements.  Finally, the third composition is $(n)$.  The set $X^{(n)}$ consists of a single, one-rowed tabloid which we will use to index the entire set $\{ 1, \dots, n \}$.  

As an example of how we will use $M^{(1,\dots,1)}$ and $M^{(k,n-k)}$ in this paper, in the next section we will consider the voting situation in which voters have been asked to fully rank a set of $n$ candidates.  We will then use that information, which we will view as an element of $M^{(1,\dots,1)}$, to assign points to each of the individual candidates.  The result may therefore be viewed as an element of $M^{(1,n-1)}$.  Furthermore, we will do all of this using a linear transformation from $M^{(1,\dots,1)}$ to $M^{(1,n-1)}$.

Now it is natural in voting to insist that a voting method not depend on how we label the candidates.  Fortunately, we may easily express this property using a group action.  More specifically, note that the symmetric group $S_n$ acts naturally on $X^\lambda$ where if $\sigma \in S_n$ and $x \in X^\lambda$, then $\sigma \cdot x$ is the tabloid one gets by applying the permutation $\sigma$ to each entry of $x$.  

The action of $S_n$ on $X^\lambda$ extends to an action of $S_n$ on $M^\lambda$ where if $\sigma \in S_n$, $\mathbf{f} \in M^\lambda$, and $x \in X^\lambda$, then 
\[
(\sigma \cdot \mathbf{f})(x) = \mathbf{f}(\sigma^{-1} \cdot x). 
\]
In other words, $M^{\lambda}$ may be viewed as a module over the group algebra $\mathbb{R}S_n$.  For example, $M^{(n)}$ corresponds to the trivial $\mathbb{R}S_n$-module.  This is because $X^{(n)}$ consists of a single tabloid, and every permutation in $S_n$ fixes this tabloid.  

In the voting example above, we can view the action of $S_n$ on $M^{(1,\dots,1)}$ and $M^{(1,n-1)}$ as the result of changing the labels on the candidates.  When the outcome of a voting method does not depend on the labels that have been assigned to the candidates, we say that the voting method is \emph{neutral} (see, for example,  \cite{TaylorManip,SaariBGOV}, among many other references for general voting theory concepts).  Insisting that a voting method based on a linear transformation $T: M^{(1,\dots,1)} \to M^{(1,n-1)}$ be neutral then simply becomes the requirement that $T$ must be an $\mathbb{R}S_n$-module homomorphism.  We will focus only on voting methods that are neutral in this paper.  

Fortunately, the representation theory of $M^\lambda$ is well-understood.  The irreducible $\mathbb{R}S_n$-modules are parametrized by the partitions of $n$, and we denote the irreducible module corresponding to the partition $\mu$ by $S^\mu$.  These are the well-known \emph{Specht modules} where, for example, $S^{(n)}$ corresponds to the one-dimensional trivial $\mathbb{R}S_n$-module.  

If $\lambda$ is a partition of $n$ (and not just a composition), then $M^\lambda$ is isomorphic to a direct sum of Specht modules
\[
M^\lambda \cong \bigoplus_\mu  \kappa_{\mu \lambda}S^\mu
\]
where the $\kappa_{\mu \lambda}$ are \emph{Kostka numbers} and are used here record the multiplicity of each Specht module in $M^\lambda$ (see, for example, Sections 2.3 and 2.11 of~\cite{Sagan}).

Note that if $\lambda$ is a composition of $n$ and not necessarily a partition, then the module $M^\lambda$ is easily seen to be isomorphic to $M^{\overline{\lambda}}$ by simply reordering the rows of the tabloids in $X^\lambda$.  We may therefore just as easily work with compositions as with partitions of $n$ when dealing with real-valued functions (e.g., voting data) defined on sets (e.g., rankings of candidates) that are indexed by tabloids.  

The modules $M^{(1,\dots,1)}$, $M^{(k,n-k)}$, and $M^{(n)}$ have particularly straightforward decompositions in terms of Specht modules.  For example, $M^{(1,\dots,1)}$ is isomorphic to the regular $\mathbb{R}S_n$-module, and thus
\[
M^{(1,\dots,1)} \cong \bigoplus_\mu (\dim S^\mu) S^\mu.
\]
If $1 \le k \le n/2$ (so that $(n-k,k)$ is a partition of $n$), then 
\[
M^{(k,n-k)} \cong M^{(n-k, k)} \cong S^{(n)} \oplus S^{(n-1,1)} \oplus S^{(n-2,2)} \oplus \cdots \oplus S^{(n-k,k)}.
\]
Lastly, because $M^{(n)}$ corresponds to the trivial $\mathbb{R}S_n$-module, we have that $M^{(n)} \cong S^{(n)}$ (see Chapter 7 of \cite{Diaconis} or Section 2.11 of~\cite{Sagan}).

Before leaving this section, we introduce some specialized notation and terminology that we will use later in the paper.  

First, note that $M^{(1,n-1)}$ is a direct sum of two irreducible submodules.  These submodules are easy to describe (see also p.~39 in~\cite{Sagan}).  One of the submodules, which we will denote by $U_0$, is isomorphic to $S^{(n)}$ and consists of all of the constant functions in $M^{(1,n-1)}$:
\[
U_0 = \left\{ \mathbf{f} \in M^{(1,n-1)} \ | \ \mathbf{f}(x) = \mathbf{f}(y) \textnormal{ for all } x,y \in X^{(1,n-1)} \right\}.  
\]
The other submodule, which we will denote by $U_1$, is isomorphic to $S^{(n-1,1)}$.  It is the orthogonal complement of $U_0$, and it consists of those functions whose values sum to zero:
\[
U_1 = \Set{ \mathbf{f} \in M^{(1,n-1)} \ | \ \sum_{x \in X^{(1,n-1)}} \mathbf{f}(x) = 0 }.
\]
If $\mathbf{f} \in M^{(1,n-1)}$, then we will denote its projection into $U_1$ by $\widehat{\mathbf{f}}$.  Thus $\mathbf{f} = (\mathbf{f} - \widehat{\mathbf{f}}) + \widehat{\mathbf{f}}$ where $\mathbf{f} - \widehat{\mathbf{f}} \in U_0$ and $\widehat{\mathbf{f}} \in U_1$.  In fact, if we let $\mathbf{1} \in M^{(1,n-1)}$ denote the function that assigns the value 1 to every tabloid, and let $\mathbf{1}_\mathbf{f}$ denote the projection of $\mathbf{f}$ onto $\mathbf{1}$, then $\mathbf{f} - \widehat{\mathbf{f}} = \mathbf{1}_\mathbf{f}$.  

Next, note that the action of $S_n$ on $X^\lambda$ is transitive.  In other words, given two tabloids $x, y \in X^\lambda$, there exists at least one permutation $\sigma \in S_n$ such that $\sigma \cdot x = y$.  In this case, we have that $\sigma \cdot \mathbf{f}_x = \mathbf{f}_y$.  It follows that every module homomorphism defined on $M^\lambda$ is determined by the image of any of the indicator functions $\mathbf{f}_x$.  For convenience, when we construct such maps, we will focus on the image of $\mathbf{f}_{x_0}$, where we use $x_0$ to denote the tabloid in $X^\lambda$ that appears first when the tabloids in $X^\lambda$ are listed lexicographically.  More specifically, $x_0 \in X^\lambda$ is the tabloid that contains the tableau whose entries, when read from left to right and top to bottom, are the numbers $1,\dots,n$ in that order.  

Finally, suppose $\lambda$ is a composition of $n$, and that $T$ is a module homomorphism defined on $M^\lambda$.  The module $M^\lambda$ can be written as the direct sum
\[
M^\lambda = \ker T \oplus (\ker T)^\perp
\]
of the kernel of $T$ and its orthogonal complement.  We will refer to the submodule $(\ker T)^\perp$ as the \emph{effective space} of $T$, and we will denote it by $E(T)$.  Note that the effective space $E(T)$ is isomorphic to (and also determines) the image of $T$.


\section{Voting Theory}

For algebraists interested in learning about the mathematics of voting, Donald Saari's papers \cite{SaariChaos, SaariAllThree, SaariStruct1, SaariStruct2} and books \cite{SaariGOV, SaariBGOV, SaariChaoticElections, SaariDecisions} are a fruitful place to start.  His work provides a friendly but mathematically sophisticated gateway to the subject, and although his approach is primarily geometric, linear algebra and symmetry are used to great effect.  

Some of Saari's key results rely on decomposing vector spaces of voting data into a handful of simple but meaningful subspaces.  Despite the geometric flavor of many of his papers, algebraists will easily recognize many of these subspaces as submodules of $M^{(1,\dots,1)}$.  This realization has been put to use in papers such as \cite{CrismanPermuta} and \cite{OrrisonSymmetry}, and what follows in this section can be partly viewed as an introduction to those two papers.

We begin by considering an election in which there are $n$ candidates labeled $1, \dots, n$.  Suppose each voter has been asked to rank the candidates by choosing a tabloid in $X^{(1,\dots,1)}$, where the candidate in the top row of the tabloid is their favorite, the candidate in the second row is their second favorite, and so on.  

For example, if there are $n = 3$ candidates, then each voter is being asked to choose one of the following tabloids from $X^{(1,1,1)}$:
\begin{center}
\begin{tabular}{| c |} \hline
1 \\ \hline
2 \\ \hline
3 \\ \hline
\end{tabular} \hspace{.05in},\hspace{.05in}
\begin{tabular}{| c |} \hline
1 \\ \hline
3 \\ \hline
2 \\ \hline
\end{tabular} \hspace{.05in},\hspace{.05in}
\begin{tabular}{| c |} \hline
2 \\ \hline
1 \\ \hline
3 \\ \hline
\end{tabular} \hspace{.05in},\hspace{.05in}
\begin{tabular}{| c |} \hline
2 \\ \hline
3 \\ \hline
1 \\ \hline
\end{tabular} \hspace{.05in},\hspace{.05in}
\begin{tabular}{| c |} \hline
3 \\ \hline
1 \\ \hline
2 \\ \hline
\end{tabular} \hspace{.05in},\hspace{.05in}
\begin{tabular}{| c |} \hline
3 \\ \hline
2 \\ \hline
1 \\ \hline
\end{tabular} \ .
\end{center}
In this case, a voter who chooses the third tabloid above is saying that she prefers candidate $2$ to candidate $1$ to candidate $3$.  

Let $\mathbf{p} \in M^{(1,\dots,1)}$ be the function with the property that $\mathbf{p}(x)$ is the number of voters who chose the tabloid $x$.  The function $\mathbf{p}$ is called a \emph{profile}.  Next, let $\mathbf{w} = [w_1,\dots,w_n]^t$ be a column vector in $\mathbb{R}^n$ such that $w_1 \ge \cdots \ge w_n$.  The vector $\mathbf{w}$ is called a \emph{weighting vector}, and we will use this vector of ``weights" to create a voting procedure as follows.  

For each tabloid $x \in X^{(1,\dots,1)}$, if candidate $i$ is in row $j$ of $x$, then she will be given $\mathbf{p}(x) w_j$ points.  We then sum over all of the tabloids in $X^{(1,\dots,1)}$ to find the total number of points received by candidate $i$.  The candidates who receive the most points (there might be ties) are then declared to be winners.  We refer to this voting procedure as the \emph{positional voting} procedure associated to $\mathbf{w}$.  Note that the inequality condition above ensures that a voter's ``more favored'' candidates receive at least as many points from her as her ``less favored" candidates.  

Note that $\mathbf{w} = [1, 0, \dots, 0]^t$ corresponds to the \emph{plurality} voting procedure, where voters are essentially being asked to vote for only one candidate.  When $\mathbf{w} = [1,\dots,1,0]^t$, we get the $\emph{anti-plurality}$ voting procedure, where voters are essentially being asked to vote against their least favorite candidate.  The \emph{Borda count}, where a candidate receives $(i-1)$ points for each time she is ranked in the $i$th position by a voter, is given by the vector $\mathbf{w} = [n-1, n-2, \dots, 2,1,0]^t$.  

Let $\mathbf{w}^{(i)} \in M^{(1,\dots,1)}$ be defined by setting $\mathbf{w}^{(i)}(x) = w_j$ if candidate $i$ is contained in the $j$th row of the tabloid $x$.  We can use the function $\mathbf{w}^{(i)}$ to find the total number of points candidate $i$ will receive, which is
\[
\sum_{x \in X^{(1,\dots,1)}} \mathbf{p}(x) \mathbf{w}^{(i)}(x) = \langle \mathbf{p}, \mathbf{w}^{(i)} \rangle.  
\]
Furthermore, we can use the tabloids in $X^{(1,n-1)}$ to index the candidates by associating $y \in X^{(1,n-1)}$ with the single candidate in the top row of $y$.  With that in mind, let $\mathbf{c}_i \in M^{(1,n-1)}$ denote the indicator function of the tabloid corresponding to candidate $i$.  We can then create the module homomorphism $T_\mathbf{w} : M^{(1,\dots,1)} \to M^{(1,n-1)}$ by setting
\[
T_\mathbf{w} (\mathbf{f}) = \langle \mathbf{f}, \mathbf{w}^{(1)} \rangle \mathbf{c}_1 + \cdots + \langle \mathbf{f}, \mathbf{w}^{(n)} \rangle \mathbf{c}_n.
\]
The coefficient in front of $\mathbf{c}_i$ in $T_\mathbf{w}(\mathbf{p})$ is therefore simply the number of points that candidate $i$ received when we apply the positional voting procedure associated with $\mathbf{w}$ to the profile $\mathbf{p}$.  

To see the motivation behind the notation defined above, consider encoding the calculation of $T_\mathbf{w}(\mathbf{p})$ as a matrix-vector multiplication.  For example, let $n = 3$, and suppose $\mathbf{w} = [1, s, 0]^t$ for some $s$ such that $1 \ge s \ge 0$.  If the coordinate vector of $\mathbf{p}$ with respect to the usual basis of $M^{(1,1,1)}$ is $[3, 2, 0, 2, 0, 4]^t$, then the coordinate vector of $T_\mathbf{w}(\mathbf{p})$ with respect to the usual basis of $M^{(1,2)}$ is given by
\begin{equation*}
\begin{bmatrix}
1 & 1 & s & 0 & s & 0 \\
s & 0 & 1 & 1 & 0 & s \\
0 & s & 0 & s & 1 & 1
\end{bmatrix}
\begin{bmatrix}
3 \\ 2 \\ 4 \\ 2 \\ 0 \\ 3
\end{bmatrix}
=
\begin{bmatrix}
5 + 4s \\ 6 + 6s \\ 3 + 4s
\end{bmatrix}.
\end{equation*}
In this case, the functions $\mathbf{w}^{(1)}, \mathbf{w}^{(2)}, \mathbf{w}^{(3)}$ correspond to the three rows of the 3-by-6 matrix, and we are simply dotting each of these rows with the coordinate vector of $\mathbf{p}$ to find the coordinate vector $[5+4s, 6+6s, 3+4s]^t$ of points received by each of the candidates.  (Note that $T_\mathbf{w}(\mathbf{p}) = (5+4s) \mathbf{c}_1 + (6+6s)\mathbf{c}_2 + (3+4s)\mathbf{c}_3$.)  

There is, however, more to see here.  Notice that each column of the 3-by-6 matrix can be obtained by permuting the entries of the weighting vector $\mathbf{w} = [1, s, 0]^t$.  In fact, if $x_0 \in X^{(1,1,1)}$ is the tabloid with candidate $i$ in row $i$, and $\sigma \in S_n$, then the column of the matrix corresponding to the tabloid $\sigma \cdot x_0$ corresponds to the function $(1)\mathbf{c}_{\sigma(1)} + (s)\mathbf{c}_{\sigma(2)} + (0)\mathbf{c}_{\sigma(3)}$ in $M^{(1,2)}$. 

From this point on we will slightly abuse our notation by identifying the weighting vector $\mathbf{w}$ with the function $w_1 \mathbf{c}_1 + \cdots + w_n \mathbf{c}_n \in M^{(1,n-1)}$.  If $x_0 \in X^{(1,\dots,1)}$ is the tabloid with candidate $i$ in row $i$, and $\widetilde{\mathbf{f}} \in \mathbb{R}S_n$ is defined by setting $\widetilde{\mathbf{f}}(\sigma) = \mathbf{f}(\sigma \cdot x_0)$, we then have that 
\begin{align*}
T_\mathbf{w}(\mathbf{f}) 
&= \sum_{\sigma \in S_n} \mathbf{f}(\sigma \cdot x_0) (w_1 \mathbf{c}_{\sigma(1)} + \cdots + w_n \mathbf{c}_{\sigma(n)}) \\
&= \sum_{\sigma \in S_n} \widetilde{f}(\sigma) (\sigma \cdot \mathbf{w}) \\ 
&= \widetilde{\mathbf{f}} \cdot \mathbf{w}.  
\end{align*}
In other words, we may interpret the function $T_\mathbf{w}(\mathbf{f})$ as the result of the group algebra element $\widetilde{\mathbf{f}} \in \mathbb{R}S_n$ acting on the function $\mathbf{w} \in M^{(1,n-1)}$.  Finding functions $\mathbf{f} \in M^{(1,\dots,1)}$ so that $T_\mathbf{w}(\mathbf{f})$ has certain desirable properties then becomes a question of finding appropriate elements of the group algebra $\mathbb{R}S_n$ to act on $\mathbf{w} \in M^{(1,n-1)}$.  This was one of the key insights in \cite{OrrisonSymmetry}, and it leads to theorems like the following, which is essentially a special case of Theorem~1 in that paper:  

\begin{theorem}\label{linindep}
Let $n \ge 2$, and suppose $\mathbf{w}_1, \dots, \mathbf{w}_k \in U_1 \subset M^{(1,n-1)}$ form a linearly independent set of weighting vectors.  If $\mathbf{r}_1, \dots, \mathbf{r}_k \in U_1$, then there exist infinitely many functions $\mathbf{f} \in M^{(1,\dots,1)}$ such that $T_{\mathbf{w}_i}(\mathbf{f}) = \mathbf{r}_i$ for all $i$ such that $1 \le i \le k$.  
\end{theorem}

To begin to appreciate the connection between Theorem~\ref{linindep} and voting theory, it is first helpful to realize that if $\mathbf{w} = \mathbf{1}_\mathbf{w} + \widehat{\mathbf{w}}$, where $\mathbf{1}_\mathbf{w} \in U_0$ and $\widehat{\mathbf{w}} \in U_1$, then $T_\mathbf{w}(\mathbf{f}) = \widetilde{\mathbf{f}} \cdot \mathbf{w} = \widetilde{\mathbf{f}} \cdot \mathbf{1}_\mathbf{w} + \widetilde{\mathbf{f}} \cdot \widehat{\mathbf{w}}$.  Furthermore, the ordinal ranking (i.e., who comes in first, who comes in second, and so on) provided by $T_\mathbf{w}$ is going to depend only on $\widehat{\mathbf{w}}$.  After all, $\widetilde{\mathbf{f}} \cdot \mathbf{1}_\mathbf{w}$ is contained in $U_0$, and is therefore a constant function.  

Building on this insight, we will say that two weighting vectors $\mathbf{w}, \mathbf{w'} \in M^{(1,n-1)}$ are equivalent, and write $\mathbf{w} \sim \mathbf{w'}$, if and only if there exist $\alpha, \alpha' \in \mathbb{R}$ such that $\alpha > 0$ and $\mathbf{w'} = \alpha \mathbf{w} + \alpha' \mathbf{1}$.   It should be clear that in this case, $\mathbf{w'}$ and $\mathbf{w}$ will yield the same ordinal rankings.  (Also, note that if $\alpha < 0$ then the ordinal rankings given by $\mathbf{w'}$ would be the reverse of the ordinal rankings given by $\mathbf{w}$.)
The following theorem, which follows from Theorem~\ref{linindep} and Theorem 2.3.1 in \cite{SaariGOV}, highlights the usefulness of this equivalence relation on weighting vectors when dealing with positional voting procedures:  

\begin{theorem}
Let $\mathbf{w}, \mathbf{w'} \in M^{(1,n-1)}$ be weighting vectors.  The ordinal rankings of $T_\mathbf{w}(\mathbf{p})$ and $T_\mathbf{w'}(\mathbf{p})$ will be the same for all profiles $\mathbf{p} \in M^{(1,\dots,1)}$ if and only if $\mathbf{w} \sim \mathbf{w'}$.  
\end{theorem}

In other words, if the weighting vectors $\mathbf{w}$ and $\mathbf{w'}$ are not equivalent, then the ordinal rankings of the outcomes of elections that use $\mathbf{w}$ and $\mathbf{w'}$ might differ, even if the same profile $\mathbf{p} \in M^{(1,\dots,1)}$ is used in both elections.  In fact, by Theorem~\ref{linindep}, the results of such elections can be arbitrarily different from one another.  This follows from the fact that if $\mathbf{w}$ and $\mathbf{w'}$ are not equivalent, then their projections into $U_1$ will be linearly independent.

Another way of appreciating the effect that Theorem~\ref{linindep} can have on our understanding of positional voting procedures is to realize that the effective spaces of $T_\mathbf{w}$ and $T_\mathbf{w'}$ share very little in common unless $\mathbf{w} \sim \mathbf{w'}$.  This is because if $\mathbf{w} \in U_1$ is nonzero, then $E(T_\mathbf{w}) \cong S^{(n-1,1)}$, and thus $E(T_\mathbf{w})$ is an irreducible submodule of $M^{(1,\dots,1)}$.  This allows us to easily prove theorems like the following, which is Theorem~4 in \cite{OrrisonSymmetry}:

\begin{theorem}\label{effective}
Let $\mathbf{w}, \mathbf{w'} \in U_1 \subset M^{(1,n-1)}$ be nonzero weighting vectors.  Then $E(T_\mathbf{w}) = E(T_\mathbf{w'})$ if and only if $\mathbf{w} \sim \mathbf{w'}$.  Furthermore, if $E(T_\mathbf{w}) \ne E(T_\mathbf{w'})$, then $E(T_\mathbf{w}) \cap E(T_\mathbf{w'}) = \{ \mathbf{0} \}$.  
\end{theorem}

In other words, if $\mathbf{w}$ and $\mathbf{w'}$ are not equivalent, then $T_\mathbf{w}$ and $T_\mathbf{w'}$ will use very different subspaces of $M^{(1,\dots,1)}$ from which to pull the information necessary to determine the outcome of an election.  The projections of a profile $\mathbf{p} \in M^{(1,\dots,1)}$ into those subspaces determine the cardinal rankings of the candidates.  We should therefore not expect any relationship between the associated ordinal rankings of the candidates.  

There are, of course, other voting methods besides positional voting procedures.  For example, a \emph{simple ranking scoring function}, or SRSF,  takes a profile $M^{(1,\dots,1)}$ and uses it to assign points to full rankings instead of individual candidates (see \cite{ConitzerRognlieXia}, and the related paper \cite{Zwicker} where these are posited as a subset of a much larger class of generalized scoring functions).  A winning \emph{ranking} (as opposed to simply a winner, which is usually how one interprets normal positional rules) is a ranking that receives at least as many points as all of the other rankings. 

As an example of an SRSF, let $x_0 \in M^{(1,\dots,1)}$ be the tabloid (i.e., full ranking) that has candidate $i$ in row $i$.  Let $\mathbf{z} \in M^{(1,\dots,1)}$ be a fixed function, and define $T_\mathbf{z}: M^{(1,\dots,1)} \to M^{(1,\dots,1)}$ by setting
\[
T_\mathbf{z}(\mathbf{f}) = \sum_{\sigma \in S_n} \mathbf{f}(\sigma \cdot x_0) (\sigma \cdot \mathbf{z}) = \widetilde{\mathbf{f}} \cdot \mathbf{z} 
\]
for all $\mathbf{f} \in M^{(1,\dots,1)}$.  Thus, if $\mathbf{p} \in M^{(1,\dots,1)}$ is a profile, then a winning ranking in our election would correspond to a tabloid $x \in M^{(1,\dots,1)}$ with the property that $T_\mathbf{z}(\mathbf{p})(x) \ge T_\mathbf{z}(\mathbf{p})(y)$ for all $y \in M^{(1,\dots,1)}$.  

How might we find a sensible function $\mathbf{z} \in M^{(1,\dots,1)}$ to create an SRSF like the one above? One way is to use a metric defined on the tabloids in $X^{(1,\dots,1)}$.  For example, if $d: X^{(1,\dots,1)} \times X^{(1,\dots,1)}$ is a metric with maximum distance $d_\textnormal{max}$, then we can set $\mathbf{z}(x) = d_\textnormal{max} - d(x,x_0)$.  A winning ranking with respect to $T_\mathbf{z}$ would then correspond to a tabloid that is ``closest" to the entire multiset of tabloids chosen by the voters (see \cite{Zwicker} for a thorough discussion of ``proximity rules").  

One of the most popular metrics on $X^{(1,\dots,1)}$ is the \emph{Kendall tau distance}, which measures the number of pairwise disagreements between two rankings. To explain, let $i, j \in \{ 1, \dots, n \}$ where $i \ne j$.  Let $\mathbf{a}_{ij} \in M^{(1,\dots,1)}$ be the function defined by setting $\mathbf{a}_{ij}(x) = 1$
whenever candidate $i$ is ranked above $j$ in $x$ (i.e., in the tabloid $x$, $i$ is in a higher row than $j$), and setting $\mathbf{a}_{ij}(x) = 0$ otherwise.  Kendall's tau distance is then given by 
\[
d(x,y) = \binom{n}{2} - \sum_{i \ne j} \mathbf{a}_{ij}(x) \mathbf{a}_{ij}(y).  
\]
The maximum distance between rankings in this case is $\binom{n}{2}$, and therefore we could define the function $\mathbf{z} \in M^{(1,\dots,1)}$ by setting $\mathbf{z}(x) = \sum_{i \ne j} \mathbf{a}_{ij}(x) \mathbf{a}_{ij}(x_0)$.  The resulting voting procedure given by $T_\mathbf{z}$ then becomes the well-known and well-studied \emph{Kemeny rule} \cite{Kemeny}.  To distinguish this case, we will denote this particular instance of $T_\mathbf{z}$ by $K:  M^{(1,\dots,1)} \to M^{(1,\dots,1)}$.  

Perhaps not surprisingly, we can use the $\mathbf{a}_{ij}$ to apply $K$ to a profile $\mathbf{p} \in M^{(1,\dots,1)}$.  More specifically, with a little bit of effort, and using the fact that the Kendall tau distance is invariant under the action of $S_n$, one can show that
\[
K(\mathbf{p}) = \sum_{i \ne j} \langle \mathbf{p}, \mathbf{a}_{ij} \rangle \mathbf{a}_{ij}
\]
(see Proposition~1 in \cite{Sibony} for a similar expression). On the other hand, it might not be clear how a voting procedure like the Kemeny rule is related to positional voting procedures.  We describe such a relationship next.  

First, it turns out that we can turn any positional voting procedure into an SRSF as follows.  Let $\mathbf{w} \in M^{(1,n-1)}$ be a weighting vector, and let $\mathbf{b} \in M^{(1,n-1)}$ be the weighting vector for the Borda Count: 
\[
\mathbf{b} = (n-1) \mathbf{c}_1 + (n-2) \mathbf{c}_2 + \dots + (1)\mathbf{c}_{n-1} + (0)\mathbf{c}_n.  
\]
Note that although we will use $\mathbf{b}$ in the construction below, any weighting vector whose weights are strictly decreasing would also work.  Next, recall that 
\[
T_\mathbf{w} (\mathbf{f}) = \langle \mathbf{f}, \mathbf{w}^{(1)} \rangle \mathbf{c}_1 + \cdots + \langle \mathbf{f}, \mathbf{w}^{(n)} \rangle \mathbf{c}_n.
\]
If we compose $T_\mathbf{w}$ with the adjoint $T_\mathbf{b}^*: M^{(1,n-1)} \to M^{(1,\dots,1)}$, then we create the map $T_\mathbf{b}^* \circ T_\mathbf{w}: M^{(1,\dots,1)} \to M^{(1,\dots,1)}$ where
\[
(T_\mathbf{b}^* \circ T_\mathbf{w}) (\mathbf{f}) = 
\langle \mathbf{f}, \mathbf{w}^{(1)} \rangle \mathbf{b}^{(1)} + \cdots + \langle \mathbf{f}, \mathbf{w}^{(n)} \rangle \mathbf{b}^{(n)}.
\]
It turns out that if we view $T_\mathbf{b}^* \circ T_\mathbf{w}$ as an SRSF, then it will return the ranking that agrees with the ranking given by $T_\mathbf{w}$ (see Proposition 1 in \cite{ConitzerRognlieXia}).  

For such procedures, individual candidates may then be ranked based on the order they appear in a winning ranking, recovering the winning candidate(s) in the usual interpretation of $T_\mathbf{w}$.  However, for more general SRSFs (including the Kemeny rule), this may not always be meaningful.  For example, suppose there are three candidates $A$, $B$, and $C$, and we are given a profile where two voters prefer $ABC$, two voters prefer $CAB$, and only one voter prefers $BCA$.  The Kemeny rule then produces a tie between the rankings $ABC$ and $CAB$.

Next, let $\mathbf{c}_{ij} \in M^{(1,1,n-2)}$ be the indicator function corresponding to the tabloid that has $i$ in the top row and $j$ in the second row.  We will call the map $P: M^{(1,\dots,1)} \to M^{(1,1,n-2)}$ given by
\[
P(\mathbf{f}) = \sum_{i \ne j} \langle \mathbf{f}, \mathbf{a}_{ij} \rangle \mathbf{c}_{ij} 
\]
the \emph{pairs map} defined on $M^{(1,\dots,1)}$.  Note that if $\mathbf{p} \in M^{(1,\dots,1)}$ is a profile, then $P(\mathbf{p})$ simply catalogues the number of times each candidate was ranked over another candidate.  

For example, when $n = 3$, $M^{(1,\dots,1)} = M^{(1,1,n-2)}$, and with respect to the usual basis of $M^{(1,1,1)}$, the pairs map $P$ can be encoded as the matrix
\[
[P] = 
\begin{bmatrix}
1 & 1 & 0 & 0 & 1 & 0 \\
1 & 1 & 1 & 0 & 0 & 0 \\
0 & 0 & 1 & 1 & 0 & 1 \\
1 & 0 & 1 & 1 & 0 & 0 \\
0 & 0 & 0 & 1 & 1 & 1 \\
0 & 1 & 0 & 0 & 1 & 1 \\
\end{bmatrix}.
\]
Note, for example, that the first row of this matrix corresponds to $\mathbf{a}_{12}$, the third row corresponds to $\mathbf{a}_{21}$, and the sum of these two rows corresponds to the constant all-ones function in $M^{(1,1,1)}$.  

The pairs map $P$ is not surjective.  In particular, the codomain has the decomposition
\[
M^{(1,1,n-2)} \cong S^{(n)} \oplus 2 S^{(n-1)} \oplus S^{(n-2,2)} \oplus S^{(n-2,1,1)}
\]
but the effective space of $P$ is isomorphic to $S^{(n)} \oplus S^{(n-1,1)} \oplus S^{(n-2,1,1)}$ (see p.~682 in \cite{OrrisonSymmetry}).  Let $W_0$, $W_1$, and $W_2$ denote the corresponding subspaces in $M^{(1,\dots,1)}$, where $W_0 \cong S^{(n)}$, $W_1 \cong S^{(n-1,1)}$, $W_2 \cong S^{(n-2,1,1)}$, and 
\[
E(P) = W_0 \oplus W_1 \oplus W_2.  
\]
Note that if $P^*: M^{(1,1,n-2)} \to M^{(1,\dots,1)}$ is the adjoint of $P$, then 
\[
(P^* \circ P) (\mathbf{f}) = \sum_{i \ne j} \langle \mathbf{f}, \mathbf{a}_{ij} \rangle \mathbf{a}_{ij}.
\]
In other words, $K = P^* \circ P$.  This should be compared with $T_\mathbf{b}^* \circ T_\mathbf{w}$ defined above, particularly in the case when $\mathbf{w}=\mathbf{b}$ corresponds to the Borda count.

If necessary, we can create the matrix encoding of $K$ with respect to the usual basis of $M^{(1,\dots,1)}$ by simply taking the product of the matrix encodings of $P^*$ and $P$.  For example, when $n = 3$, the matrix encoding of $K$ is 
\[
[K] = 
\begin{bmatrix}
3 & 2 & 2 & 1 & 1 & 0 \\
2 & 3 & 1 & 0 & 2 & 1 \\
2 & 1 & 3 & 2 & 0 & 1 \\
1 & 0 & 2 & 3 & 1 & 2 \\
1 & 2 & 0 & 1 & 3 & 2 \\
0 & 1 & 1 & 2 & 2 & 3 \\
\end{bmatrix}
\]
which is the product $[P]^t [P]$.  Notice the columns of $[K]$ are all permutations of the first column of $[K]$, which corresponds to the function $\mathbf{z} \in M^{(1,1,1)}$ that is based (as was described above) on the Kendall tau distance.  

Interestingly, it turns out that the effective space $E(T_\mathbf{b})$ of the Borda Count is $W_0 \oplus W_1$ (see, for example, Section~7 of~\cite{OrrisonSymmetry}).  This is perhaps much more understandable once we realize that 
\[
\mathbf{b}^{(i)} = \sum_{j: \ i \ne j} \mathbf{a}_{ij}.
\]
Furthermore, by Theorem~\ref{effective}, the only nontrivial positional voting procedures that contain $W_1$ in their effective spaces are those whose weighting vectors are equivalent to $\mathbf{b}$.  This begins to explain theorems like the following, which is a combination of parts of Theorem~3 and Theorem~4 in \cite{SaariMerlin1}:

\begin{theorem}
For $n \ge 3$ candidates, the Borda count always ranks the Kemeny rule top-ranked candidate strictly above the Kemeny rule bottom ranked-candidate.  Conversely, the Kemeny rule ranks the Borda count top-ranked candidate strictly above the Borda count bottom-ranked candidate.  For any positional voting method other than the Borda count, however, there is no relationship between the Kemeny rule ranking and the positional ranking.  
\end{theorem}

Using the pairs map $P$, we can say more about the relationship between the Borda count and the Kemeny rule.  More specifically, because $K = P^* \circ P$ is self-adjoint, we know it is orthogonally diagonalizable.  In fact, the eigenspaces of $K$ are $W_0$, $W_1$, $W_2$, and $(W_0 \oplus W_1 \oplus W_2)^\perp$ with  eigenvalues $\kappa_0 = \frac{n!}{2}\binom{n}{2}$, $\kappa_1 = \frac{(n+1)!}{6}$, $\kappa _2 = \frac{n!}{6}$, and $0$, respectively (see Theorem~3 in~\cite{Sibony}).  

If we let $T_i: M^{(1,\dots,1)} \to W_i$ denote the orthogonal projection onto $W_i$, we then have that
\[
K = \kappa_0 T_0 + \kappa_1 T_1 + \kappa_2 T_2.  
\]
Furthermore, if $\beta_0 = \frac{(n-1) n!}{2} \binom{n}{2}$, $\beta_1 = \frac{n (n+1)!}{12}$, and $\beta_2 = 0$ then it is possible to show that 
\[
T_\mathbf{b}^* \circ T_\mathbf{b} = \beta_0 T_0 + \beta_1 T_1+\beta_2 T_2.  
\]
The Borda count and Kemeny rule may therefore be viewed as members of the same family of SRSFs whose maps all have the form
\[
K_{(\gamma_0, \gamma_1, \gamma_2)} = \gamma_0 T_0 + \gamma_1 T_1 + \gamma_2 T_2
\]
where $\gamma_0, \gamma_1, \gamma_2 \in \mathbb{R}$.  

When using $K_{(\gamma_0, \gamma_1, \gamma_2)}$, the resulting ordering of the rankings in $X^{(1,\dots,1)}$ depends only on the ratio $\gamma_2/\gamma_1$.  After all, $\gamma_0 T_0$ only contributes scalar multiples of the all-ones function to the outcome.  This fact is exploited in \cite{CrismanPermuta} by setting $\gamma_0 = 0$ and focusing on a one-parameter family of procedures that is interpreted as being ``between" the Borda count and Kemeny rule.  Such a family allows us, for example, to better understand when properties like being susceptible to the ``no-show paradox" can arise in voting (see Proposition 5.16 in~\cite{CrismanPermuta}).  

By describing everything up to this point in terms of $M^{(1,\dots,1)}$ and $M^{(1,n-1)}$, we hope it is clear how one might extend all of these ideas to more general settings.  For example, suppose $\lambda = (\lambda_1, \dots, \lambda_m)$ is a partition of $n$.  We could ask the voters to rank the candidates by choosing a tabloid in $X^\lambda$, where the candidates in the top row of their tabloid are their top favorite $\lambda_1$ candidates, the candidates in the second row are their next favorite $\lambda_2$ candidates, and so on.  

When voters are asked to provide rankings by choosing a tabloid from $X^{(1,\dots,1)}$, we say they are giving \emph{full rankings} of the candidates.  If they are choosing tabloids from $X^\lambda$ where $\lambda \ne (1,\dots,1)$, then we say they are giving \emph{partial rankings} (of shape $\lambda$) of the candidates.  We can therefore ask if any of the theorems above can be extended to partial rankings.  For some of the theorems, the answer is yes.  In fact, the paper by Daugherty et al.~\cite{OrrisonSymmetry} focused primarily on such theorems for positional voting.  

Notice also that the maps $T_\mathbf{w}$ and $T_\mathbf{z}$ were defined in essentially the same way.  In both cases, we chose a vector $\mathbf{v} \in M^\mu$ for some partition $\mu$, and then  defined $T_\mathbf{v}: M^{(1,\dots,1)} \to M^{\mu}$ by setting $T_\mathbf{v}(\mathbf{f}) = \widetilde{f} \cdot \mathbf{v}$.  We then used $T_\mathbf{v}$ and a profile $\mathbf{p} \in M^{(1,\dots,1)}$ to find the ``best" tabloid in $M^\mu$ based on  ``points" that $T_\mathbf{v}(\mathbf{p})$ assigned to the standard basis vectors of $M^\mu$.  What would happen if we were to extend this construction to situations in which voters choose tabloids in $X^\lambda$ and points are then assigned to tabloids in $X^\mu$?  Would we learn anything new?  For readers interested in exploring this question, Section 3 in \cite{OrrisonSymmetry} might be a good place to start.  It discusses maps from $M^\lambda$ to $M^{(1,n-1)}$ and what happens to positional voting procedures when voters are asked to provide partial (instead of full) rankings of shape $\lambda$.  

Another possible direction one might take is to ask what would happen if we were to replace the pairs map in the definition of the Kemeny rule with maps that catalogue information about triples of candidates, or quadruples of candidates, and so on.  Doing so would then place the Kemeny rule in another potentially interesting family of voting procedures.  We are not aware of any work in this direction, but anyone interested in pursuing such a project would almost certainly benefit from the discussions about inversions in \cite{GrossmanMinton} and \cite{Marden}.

Finally, there are \emph{maximum likelihood estimator} procedures which give the ``most likely" ranking of the candidates, assuming that the voters all meant to choose the same ranking but their votes were corrupted by a noise model.  Young \cite{Young} shows that the Borda count is an MLE when the desired outcome is a \emph{single} winner, and that the Kemeny rule is an MLE when the desired outcome is a full ranking.  In \cite{ProcacciaMLE}, the question of which procedure is an MLE when the desired outcome is a top pair, top triple, and so forth is investigated.  Given the connection between SRSFs and MLEs found in \cite{ConitzerRognlieXia}, we suspect that an algebraic approach to MLE procedures has tremendous promise.


\section{Game Theory}

For algebraists interested in learning about (transferable utility cooperative) game theory, the pioneering work of Kleinberg and Weiss \cite{KWAlgGames, KWWeakValues, KWNullShapley, KWMembership, KWOrthogonal} will quickly make them feel welcome.  Consider, for example, the first few lines from the introduction of~\cite{KWAlgGames}:
\begin{quote}
Game theory and algebra become inextricably intertwined once one recognizes that the notion of a permutation of players gives rise to a representation of the symmetric group in the space of automorphisms of $\mathcal{G}$, the vector space of games (in characteristic function form).  

In this paper we use this observation to attempt to turn the usual approach to game-theoretic problems on its head by analyzing the space of games as an algebraist might.  In this way, we can see if the mathematical structure of $\mathcal{G}$ has anything to say about what constructs are significant from a game-theoretic point of view.  
\end{quote}

In this section, we focus on the role that tabloids play when navigating Kleinberg and Weiss's work, as well as the recent work of Hern\'andez-Lamoneda, Ju\'arez, and S\'anchez-S\'anchez~\cite{GameTheoryRepsHLJSS} and S\'anchez-P\'erez~\cite{CoopGameRepsSP}.  We then take advantage of our use of tabloids to suggest possible extensions of their work.  

We begin with some notation and terminology.  Let $n \ge 2$, and suppose we have $n$ players.  We will label the players with the numbers $1, \dots, n$, and we will let $N = \{1, \dots, n \}$ be the set of all of the players.  

A \emph{cooperative game} for $N$ is a real-valued function defined on the set of all subsets of $N$, with the convention that the empty set is always mapped to $0$.  The set of all such games 
\[
\mathcal{G} = \{ v: 2^N \to \mathbf{R} \ | \ v(\emptyset) = 0 \}
\]
models ways in which various subsets (``coalitions") can be allotted a utility (``value") based on the extent to which they are perceived to contribute to the entire set of players (``the grand coalition").  

One of the defining challenges in cooperative game theory is to find useful and meaningful \emph{solution concepts}, which are functions defined on $\mathcal{G}$ that are used to determine a ``payoff" for each individual player.  More specifically, if we index the individual players by tabloids in $X^{(1,n-1)}$ as we did in the last section on voting theory, then we can think of a solution concept as a function
\[
\varphi: \mathcal{G} \to M^{(1,n-1)}
\]
where if $v \in \mathcal{G}$ and $i \in N$, then the value that $\varphi$ associates to player $i$ is the coefficient in front of $\mathbf{c}_i$ in $\varphi(v)$. If we denote this value by $\varphi(v)_i$, then 
\[
\varphi(v) = \varphi(v)_1 \mathbf{c}_1 + \cdots + \varphi(v)_n \mathbf{c}_n. 
\]
It is interesting to note that, algebraically, a solution concept looks similar to a positional voting procedure, in that the payoff that a solution concept assigns to a player is similar to the points that a positional voting procedure assigns to a candidate.  

The vector space $\mathcal{G}$ of games is of course much more than just a vector space.  It is an $\mathbb{R}S_n$-module, where for all $\sigma  \in S_n$ and $v \in \mathcal{G}$,  
\[
(\sigma \cdot v) (S) = v(\sigma^{-1} \cdot S)
\]
for all subsets $S  \in 2^N$, and where $\sigma^{-1} \cdot S$ is the set you obtain by applying $\sigma^{-1}$ to each element of $S$.  Furthermore, the irreducible submodules of $\mathcal{G}$ are straightforward to describe.  

First, note that we can write $\mathcal{G} = \mathcal{G}_1  \oplus \mathcal{G}_2 \cdots \oplus \mathcal{G}_{n-1} \oplus \mathcal{G}_n$, where $\mathcal{G}_k$ is the subspace of all games that assign the value $0$ to any subset of players that does not have size $k$.  Next, if we index the $k$-element subsets of players using tabloids in $X^{(k,n-k)}$, then it becomes clear that 
\[
\mathcal{G}_k \cong M^{(k,n-k)}.
\]
When $k \le n/2$, we know that $M^{(k,n-k)} \cong S^{(n)} \oplus S^{(n-1, 1)} \oplus \cdots \oplus S^{(n-k,k)}$.  We also know that $M^{(k,n-k)} \cong M^{(n-k,k)}$.  Thus, if we let $U_j^k$ denote the submodule of $\mathcal{G}_k$ that is isomorphic to the Specht module $S^{(n-j,j)}$, we then have
\begin{align*}
\mathcal{G}_1 &= U_0^1 \oplus U_1^1 \\
\mathcal{G}_2 &= U_0^2 \oplus U_1^2 \oplus U_2^2 \\
\mathcal{G}_3 &= U_0^3 \oplus U_1^3 \oplus U_2^3 \oplus U_3^3 \\
& \ \ \vdots \\
\mathcal{G}_{\lfloor n/2 \rfloor} &= U_0^{\lfloor n/2 \rfloor} \oplus U_1^{\lfloor n/2 \rfloor} \oplus U_2^{\lfloor n/2 \rfloor} \oplus \cdots \oplus U_{\lfloor n/2 \rfloor}^{\lfloor n/2 \rfloor} \\
& \ \ \vdots \\
\mathcal{G}_{n-2} &= U_0^{n-2} \oplus U_1^{n-2} \oplus U_2^{n-2} \\
\mathcal{G}_{n-1} &= U_0^{n-1} \oplus U_1^{n-1} \\
\mathcal{G}_n &= U_0^n. \\
\end{align*}
We may therefore write $\mathcal{G}$ as a direct sum of the irreducible $U_j^k$, where 
\[
\mathcal{G} = \left( \bigoplus_{k = 1}^{\lfloor n/2 \rfloor} \bigoplus_{j = 0}^k U_j^k   \right) \bigoplus \left(  \bigoplus_{k > \lfloor n/2 \rfloor}^n  \bigoplus_{j=0}^{n-k} U_j^k \right)
\]
which is essentially equation (3) in~\cite{KWAlgGames}.

The above decomposition is helpful in cooperative game theory because the consensus since \cite{Shapley} is that \emph{linear} and \emph{symmetric} solution concepts are particularly important, and these solution concepts coincide with the set of all $\mathbb{R}S_n$-module homomorphisms from $\mathcal{G}$ to $M^{(1,n-1)}$.  Moreover, because the $\mathcal{G}_k$ have such simple decompositions into irreducible submodules, these homomorphisms are easy to describe.  

To explain, first recall that $M^{(1,n-1)}$ is the direct sum of $U_0$ and $U_1$, where $U_0 \cong S^{(n)}$ and $U_1 \cong S^{(n-1,1)}$.  By Schur's lemma, this means that any module homomorphism from $\mathcal{G}_k \cong M^{(k,n-k)}$ to $M^{(1,n-1)}$ must contain $U_j^k$ in its kernel for all $j \ge 2$.  We may therefore describe any homomorphism from $\mathcal{G}$ to $M^{(1,n-1)}$ by focusing only on homomorphisms defined on the $U_0^k$ and $U_1^k$ as follows.  

Let $T_0^k: U_0^k \to U_0$ and $T_1^k: U_1^k \to U_1$ be any fixed isomorphisms.  By a slight abuse of notation, we can view these isomorphisms as homomorphisms from $\mathcal{G}$ to $M^{(1,n-1)}$ by identifying them with $\iota_0 \circ T_0^k \circ p_0^k$ and $\iota_1 \circ T_1^k \circ p_1^k$, respectively, where $p_0^k$ is the projection map onto $U_0^k$, $p_1^k$ is the projection map onto $U_1^k$, and $\iota_0$ and $\iota_1$ are the usual inclusion maps into $M^{(1,n-1)}$.  

By Schur's lemma, every module homomorphism from $\mathcal{G}$ to $M^{(1,n-1)}$ may be expressed uniquely as a linear combination of the $T_0^k$ and $T_1^k$.  In other words, for every linear symmetric solution concept $\varphi:  \mathcal{G} \to M^{(1,n-1)}$, there exist scalars $c_0^1, \dots, c_0^n$ and $c_1^1, \dots, c_1^{n-1}$ such that 
\[
\varphi = (c_0^1 T_0^1 + \cdots + c_0^n T_0^n) + (c_1^1 T_1^1 + \cdots + c_1^{n-1} T_1^{n-1}). 
\]
Working with and describing linear symmetric solution concepts therefore becomes a matter of manipulating and communicating the $c_0^j$ and $c_1^j$.  Moreover, there are choices for the $T_0^k$ and $T_1^k$ that make certain properties of solution concepts easy to verify.  

For example, consider the following maps that appear in \cite{KWWeakValues, KWMembership}.  First, let
\[
A(v,k) = \binom{n}{k}^{-1}  \sum_{S : \ |S| = k} v(S)
\]
denote the average value received by the coalitions of size $k$, and let 
\[
\gamma(k) = \binom{n-2}{k-1}.
\]
Then define $T_0^k$ by setting $T_0^k(v)_i = k^{-1} A(v,k)$, and define $T_1^k$ by setting
\[
T_1^k (v)_i = \gamma(k)^{-1}  \sum_{S: \ |S| = k \textnormal{ and } i \in S} \left[ v(S) - A(v,k) \right].  
\]
(A somewhat different approach is taken in \cite{GameTheoryRepsHLJSS}, but the results are essentially scaled versions of the maps above.)  

As example of a property that can now be expressed in terms of the $c_0^j$ and $c_1^j$, consider the property of \emph{efficiency}, where a solution concept $\varphi$ is said to be \emph{efficient} if
\[
\varphi(v)_1 + \cdots + \varphi(v)_n = v(N)
\]
for all $v \in \mathcal{G}$.  It turns out that a linear symmetric solution concept will be efficient if and only if $c_0^1 = \cdots = c_0^{n-1} = 0$ and $c_0^n = 1$ (see Proposition~2 in~\cite{KWWeakValues}).  

As another example, a solution concept $\varphi$ is said to be a \emph{marginal value} if there exist scalars $m_1, \dots, m_n$ such that 
\[
\varphi(v)_i = \sum_{S: \ i \in S}  m_{|S|} (v(S) - v(S-i))
\]
for all $v \in \mathcal{G}$.  In this case, it turns out we can recover $\varphi$ by setting 
\[
c_0^k = k \left[  m_k \binom{n-1}{k-1} - m_{k+1} \binom{n-1}{k} \right]
\]
and 
\[
c_1^k = \gamma(k) \left[  m_k + m_{k+1} \right]
\]
where we define $m_{n+1} = 0$ (see Proposition~1 in~\cite{KWMembership}).  

The above approach to describing linear symmetric solution concepts becomes even more powerful when we ask about solution concepts that have more than one such property.  For example, using the insights above, we may easily verify that there must be precisely one linear symmetric solution concept that is also an efficient marginal value.  It is the well-studied \emph{Shapley value}~\cite{Shapley}, which is obtained by setting $c_0^1 = \cdots = c_0^{n-1} = 0$, $c_0^n = 1$, and $c_1^1 = \cdots = c_1^{n-1} = \frac{1}{n-1}$ (see Section~3 of \cite{KWMembership}; see also Section 4.3 of \cite{GameTheoryRepsHLJSS}).  

Similarly, \cite{GameTheoryRepsHLJSS} uses this type of analysis to prove a simple and intuitive criterion for a solution concept to be a ``self-dual" marginal value.  To explain, they define the \emph{duality operator} $*: \mathcal{G} \to \mathcal{G}$ by setting $(*v)(S) = v(N) - v(N - S)$ for all $v \in \mathcal{G}$.  A linear symmetric solution $\varphi$ is then said to be \emph{self-dual} if $\varphi(*v) = \varphi(v)$ for all $v \in \mathcal{G}$.  It turns out that a linear symmetric marginal value will be self-dual if and only if the $m_k$ defined above have the property that $m_j = m_{n-j-1}$ for all $j < n$ (see Proposition 6 of \cite{GameTheoryRepsHLJSS}).  

By viewing a linear symmetric solution concept as a module homomorphism $\varphi: \mathcal{G} \to M^{(1,n-1)}$, it is now easy to see how one might replace $M^{(1,n-1)}$ by $M^{(k,n-k)}$ to get linear symmetric solution concepts that assign payoffs to $k$-element subsets of players instead of individual players.  Such an idea was explored in~\cite{KWAlgGames}, but there is room for significant additional investigation, as the  conjectures in that paper seem to imply.  From a representation-theoretic point of view, we could easily begin by  including in the discussion above module homomorphisms defined on the $U_j^k$ when $j \ge 2$.  

In some sense, this would begin to address the question posed in some of the papers referenced above as to whether the common kernel of all linear symmetric solution concepts has \emph{useful} structure.  After all, this common kernel is the sum of all of the $U_j^k$ where $j \ge 2$, and each of these submodules would begin to contribute to generalized solution concepts that involved payoffs for pairs, triples, and so on.  Finding meaningful properties of the associated generalized solution concepts could be challenging, but doing so might ultimately prove to be illuminating, especially if they were to enhance our understanding of cherished properties (e.g., being a marginal value) that solution concepts might possess.  

The approach taken in \cite{GameTheoryRepsHLJSS} is extended in \cite{CoopGameRepsSP} to so-called \emph{games in partition function form}.  In this setting, the value of a coalition $S$ depends not just on the coalition, but on how the players \emph{not} in the coalition are themselves partitioned into coalitions.  In this way, an ``embedded coalition" is a pair $(S,Q)$ such that $Q$ is a (set) partition of $N$ and $S\in Q$.  The dimension of the resulting game space is much larger, and its decomposition into irreducible submodules has only been done for $n = 3$ and $n = 4$, which are (unsurprisingly) mostly $U_0$ and $U_1$ for those $n$.  Can a decomposition of the general case be described?  If so, can it be used to describe properties (e.g., efficiency) of the associated generalized solution concepts?  (See, for example, Corollary 4 in~\cite{CoopGameRepsSP}.)  

Finally, it is also possible to decompose \emph{non-cooperative games} \cite{JessieSaari, KalaiKalai}.  In the simplest version of this setting, players have strategies they individually decide upon, yielding varying payoffs.  As is shown in \cite{KalaiKalai}, we can decompose such games into ``fully competitive" zero-sum components and ``cooperative" components in which all of the players receive the same payoff.  The main theorem in \cite{JessieSaari} takes this further, while also giving a different interpretation to the zero-sum situation.  As an example of their results, they show that the space of $2\times 2$ games may be decomposed into a four-dimensional subspace that uniquely determines the Nash outcome, a two-dimensional subspace consisting solely of payoffs for each player which does not affect any strategic analysis, and a two-dimensional subspace governing other behavioral features.  It would be interesting to better understand the role that representation might play here.   Some first steps in this direction are being taken by Jessie~\cite{Jessie}.


\section{Conclusion}

As we have seen, the permutation representations arising from the action of the symmetric group on tabloids provide a unifying framework for understanding and extending foundational ideas in both voting theory and game theory.  There is, however, much more work that could be done, and we encourage interested readers to consider how they might use these and other ideas to contribute to our understanding of voting theory, game theory, and other mathematical behavioral sciences.


\bibliographystyle{amsplain}


\end{document}